\newif\ifarxiv
\arxivfalse

\newif\ifdraftmode
\draftmodefalse

\ifdraftmode
    \documentclass{article}
    \usepackage{arxiv}
\else
    \documentclass{article}
    \usepackage{arxiv}
\fi
\pdfoutput=1 

\usepackage[utf8]{inputenc}             
\usepackage[T1]{fontenc}                
\usepackage{lmodern}                    
\usepackage{bbm}                        

\usepackage[english]{babel}             
\usepackage{csquotes}                   

\ifarxiv
    \usepackage{arxiv}
\else
\fi

\usepackage[protrusion=true, expansion=true]{microtype} 
\usepackage[format=plain, labelfont=bf]{caption} 
\usepackage{ulem}                       
\usepackage{amsmath}                    
\usepackage{amssymb}                    
\usepackage{amsfonts}                   
\usepackage{amsthm}                     
\usepackage{amscd}                      
\usepackage{mathtools}                  
\mathtoolsset{showonlyrefs=true}        
\usepackage{latexsym}                   

\usepackage{cancel}                     
\usepackage{esint}                      
\usepackage{braket}                     
\usepackage{mathrsfs}                   
\usepackage{sansmath}

\usepackage{booktabs}                   
\usepackage{multirow}                   
\usepackage[dvipsnames]{xcolor}         
\usepackage{colortbl}                   
\usepackage{paralist}                   
\setdefaultenum{(I)}{(i)}{}{}           

\usepackage{graphicx}                   
\usepackage[rightcaption]{sidecap}      
\usepackage{wrapfig}
\usepackage{epstopdf}                   
\usepackage{tikz}                       
\usepackage[skins]{tcolorbox}
\usepackage[caption=false]{subfig}

\usepackage{pdfpages}                   
\usepackage{placeins}
\usepackage{enumitem}

\ifdraftmode
    \usepackage[english,colorinlistoftodos,textsize=tiny]{todonotes} 
    \makeatletter
    \renewcommand{\todo}[2][]{\@bsphack\@todo[#1]{\textcolor{black}{#2}}\@esphack\ignorespaces}
    \makeatother
    \usepackage[notref,notcite]{showkeys}
    
\else
    
    \usepackage[disable]{todonotes} 
\fi

\usepackage[ruled, vlined]{algorithm2e}

\usepackage{hyperref}

\usepackage[
    backend=biber,
    natbib=true,
    style=numeric-comp,
    sorting=none,
    maxnames=2,
    giveninits=true,
    url=false,
    arxiv=pdf,
    doi=false,
    eprint=true,
    isbn=false,
    date=terse,  
    block=ragged,
    maxbibnames=99
]{biblatex}
\theoremstyle{definition}
\newtheorem{definition}{Definition}

\newtheorem{remark}[definition]{Remark}

\usetikzlibrary{fit, shapes, arrows, patterns, backgrounds, fadings, calc, positioning, shadows, shadings}
\usetikzlibrary{decorations.shapes, decorations.pathreplacing, decorations.pathmorphing, decorations.markings}

\tikzset{>=latex}

\definecolor{dred}{rgb}{0.8,0,0}
\definecolor{dgreen}{rgb}{0,0.45,0}
\graphicspath{
	{figures/}
}

\tikzset{
    diagonal fill/.code 2 args={
        \pgfdeclareverticalshading[%
            tikz@axis@top,tikz@axis@middle,tikz@axis@bottom%
        ]{diagonalfill}{100bp}{%
            color(0bp)=(tikz@axis@bottom);
            color(50bp)=(tikz@axis@bottom);
            color(50bp)=(tikz@axis@middle);
            color(50bp)=(tikz@axis@top);
            color(100bp)=(tikz@axis@top)
        }
        
        \tikzset{shade, left color=#1, right color=#2, shading=diagonalfill}
    }
}

\tikzset{
    diagonal bar/.code 2 args={
        \pgfdeclareverticalshading[%
            tikz@axis@top,tikz@axis@middle,tikz@axis@bottom%
        ]{diagonalbar}{6bp}{%
            color(0bp)=(tikz@axis@bottom);
            color(50bp)=(tikz@axis@bottom);
            color(50bp)=(tikz@axis@middle);
            color(50bp)=(tikz@axis@top);
            color(1000bp)=(tikz@axis@top)
        }
        
        \tikzset{shade, left color=#1, right color=#2, shading=diagonalbar}
    }
}

\tikzset{Basic/.style={
	rectangle,
	inner sep=2pt,
	minimum width=2.2em,
	minimum height=2.2em,
	rounded corners,
	line width=0.3mm
}}

\tikzset{nix/.style={
	Basic,
    draw=none,
	fill=none
}}

\tikzset{T/.style={
    Basic,
	draw=black, 
	text=black,
	fill=black!10
}}

\tikzset{GT/.style={
    Basic,
	draw=gray,
	text=gray,
	fill=black!5
}}

\tikzset{DT/.style={
	T,
	diagonal bar={black!40}{black!10},
	shading angle=45
}}

\tikzset{GDT/.style={
	GT,
	diagonal bar={black!20}{black!5},
	shading angle=45
}}

\tikzset{LT/.style={ 
	T,
	diagonal fill={black!40}{black!10}, 
	shading angle=45 
}}

\tikzset{LTI/.style={ 
	T,
	diagonal fill={black!10}{black!40}, 
	shading angle=45 
}}

\tikzset{RT/.style={
	T, 
	diagonal fill={black!40}{black!10},
	shading angle=-45
}}

\tikzset{RTI/.style={
	T, 
	diagonal fill={black!10}{black!40},
	shading angle=-45
}}

\tikzset{GLT/.style={ 
	GT,
	diagonal fill={black!20}{black!5},
	shading angle=45 
}}

\tikzset{GRT/.style={ 
	GT,
	diagonal fill={black!20}{black!5},
	shading angle=-45 
}}

\newcounter{sarrow}

\DeclareMathOperator*{\argmin}{arg\,min}

\newcommand{\code}[1]{\texttt{#1}}
\let\oldbullet\bullet
\newlength{\raisebulletlen}
\setbox1=\hbox{$\bullet$}\setbox2=\hbox{\tiny$\bullet$}
\renewcommand\bullet{\raisebox{\raisebulletlen}{\,\tiny$\oldbullet$}\,}

\let\oldcolon\colon
\renewcommand{\colon}{\,\oldcolon\,}

\usepackage{mdframed}

\usetikzlibrary{positioning}
\tikzset{
    position/.style args={#1:#2 from #3}{
        at=(#3), anchor=#1+180, shift=(#1:#2)
    }
}

\title{Approximative Policy Iteration for Exit Time Feedback Control Problems driven by Stochastic Differential Equations using Tensor Train format}
\date{\today}

\author{
  Konstantin Fackeldey \\
            Technische Universit\"at Berlin\\
            Strasse des 17. Juni 135 \\
            10623 Berlin, Germany \\
  \texttt{fackeldey@math.tu-berlin.de} \\
   \And
  Mathias Oster \\
            Technische Universit\"at Berlin\\
            Strasse des 17. Juni 135 \\
            10623 Berlin, Germany \\
  \texttt{oster@math.tu-berlin.de} \\
   \And
  Leon Sallandt\\
            Technische Universit\"at Berlin\\
            Strasse des 17. Juni 135 \\
            10623 Berlin, Germany \\
  \texttt{sallandt@math.tu-berlin.de} \\
  \And
  Reinhold Schneider\\
            Technische Universit\"at Berlin\\
            Strasse des 17. Juni 135 \\
            10623 Berlin, Germany \\
  \texttt{schneidr@math.tu-berlin.de} \\
}

\makeatletter
\AtBeginDocument{%
  \patchcmd{\lbx@us@mkdaterangetrunc@short}
    {\endngroup}
    {\endgroup}
    {}{}%
}
\makeatother

\begin{document}
\maketitle

\begin{abstract}
    We consider a stochastic optimal exit time feedback control problem. 
    The Bellman equation is solved approximatively via the Policy Iteration algorithm on a polynomial ansatz space by a sequence of linear equations. 
    As high degree multi-polynomials are needed, the corresponding equations suffer from the curse of dimensionality even in moderate dimensions. 
    We employ tensor-train methods to account for this problem.
    The approximation process within the Policy Iteration is done via a Least-Squares ansatz and the integration is done via Monte-Carlo methods.
    Numerical evidences are given for the (multi dimensional) double well potential and a three-hole potential.
\end{abstract}

\section{Introduction}
Optimal control of ordinary differential equations (ODE) is a field of mathematics and engineering, where we minimize a cost functional constrained by a controlled ODE. 
Substituting the ODE by a stochastic differential equation (SDE) we obtain a stochastic optimal control problem. 
The inherent structure of the cost functional determines the behavior of the optimization problem. 
Within this context several formulations have been investigated. Among them are finite and infinite horizon problems and exit time problems.
For the latter, one determines the optimal control to reach a predefined exit set with respect to the cost functional.
An inherent difficulty is that the stopping time is not known in advance and depends not only on the control but also on the stochastic process.

Optimal control problems of stochastic processes have been utilized in various fields of applications, such as finance, engineering or molecular dynamics, see e.g.\cite{Fleming1999, OptimalControl-Engrg, SchuetteWinkelmannHartmann2012}. 
Subsequent to the latter in \cite{HartmannRichterSchuetteetal.2017, hartmann2019variational}
this optimal control setting has been applied to
the characterization of free energy of an uncontrolled dynamical system. Different numerical methods have been developed and are now widely used and further investigated, see e.g. \cite{Infinite1,Barron_93,Damm_17,Hartmann_2012,nsken2020solving} .
A popular approach is approximating the value function by solving either the Bellman or the Hamilton-Jacobi-Bellman (HJB) equation.

For low dimensions $ n \leq 3$ the HJB equations, corresponding to deterministic and stochastic optimal control problems, have been treated by Finite Element  \cite{smears} and Finite Difference methods \cite{Zidani_1,Zidani_2}, 
including Semi-Lagrangian methods \cite{Falcone,Falcone1987, SemiLagrangianStochastic,SemiLagranigian}. 
These methods are based on grids and are facing the curse of dimensions, which prevents the treatment of large spatial dimensions $n$. 
Popular approaches to get rid of the curse of dimensions are sparse grids \cite{Garcke.Kroener:2017,Garcke1},
Tensor trains combined with Galerkin \cite{TensorKunisch} or minimal residual \cite{oster2019approximating} methods, or using Max-Plus
algebra \cite{Akian}
, see also \cite{DataHJB, VIM, Beard} for further ideas.
Nowadays deep neural networks (DNN) have become an attractive tool for solving the HJB \cite{Han2016DeepLA,Han8505,pham2020neural}. 
In our approach we address the high-dimensionality by using Monte-Carlo integration and Tensor Train (TT) formats and the non-linearity by using the Policy Iteration algorithm \cite{howard, bellman1955, bellmann61}.



\begin{wrapfigure}{r}{9cm}
    \centering
    \includegraphics[width = 0.9\linewidth]{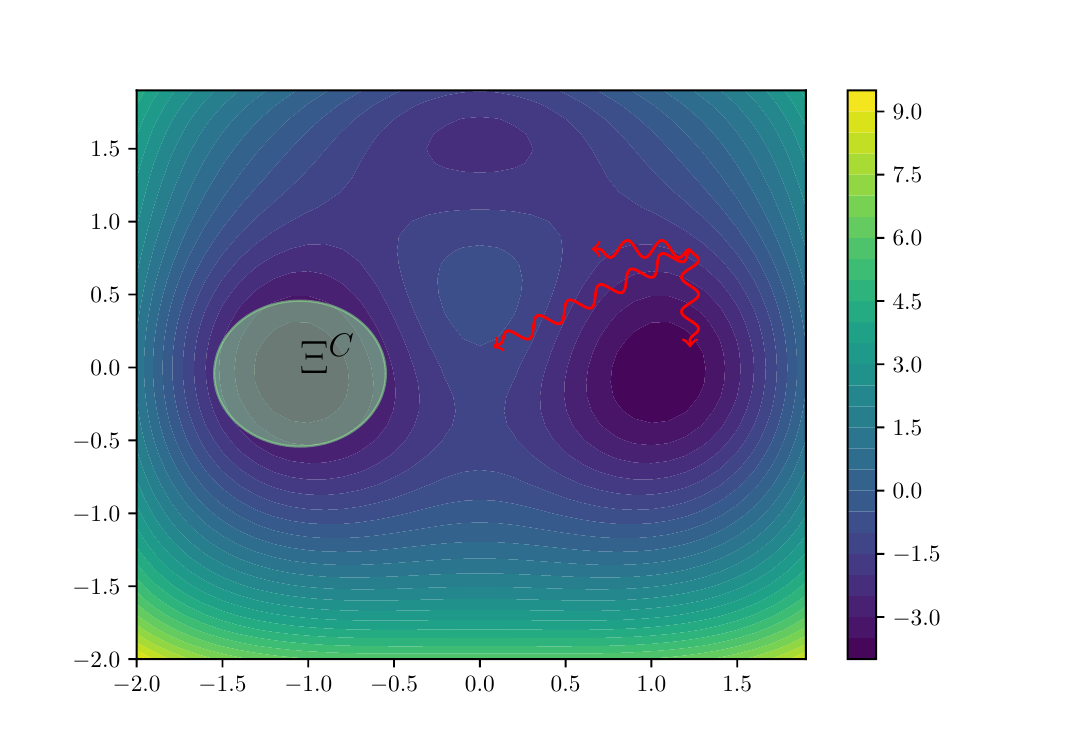}
    \caption{\label{fig:optControlProblem} At a given point on the right hand side different trajectories of the stochastic process are visualized by the red arrows. We now seek for the optimal control driving the process out of the set $\Xi$.}
    
\end{wrapfigure}
Here, we consider some stochastic process $X_t$ starting in $X_0 = x$ on a potential landscape (cp Fig \ref{fig:optControlProblem}). Our goal is to modify $X_t$ by a control $u_t$ such that the process exits a predefined set $\Xi$. Of course, there are many possible controls, driving the system out of $\Xi$. Our goal is it to find the optimal control with respect to a certain cost function $\mathcal J = \mathcal J (x, u)$.

It is worth mentioning, that the infinitesimal version of the Bellman equation is the Hamilton-Jacobi-Bellman (HJB) equation. If the structure of the stochastic part is modeled by  a Wiener process, it appears within the HJB as a Laplace operator \citep{StochContrTheo_Nisio,prato_zabczyk_1992_StochEqInfDim,Markov_Fleming_Soner, Fabbri_Gozzi_Stoch_OptConInfDim}. 
Solving the HJB could be realized by adapting the method in \cite{TensorKunisch}, where the deterministic HJB has been approximated by using Tensor Trains. Note that other function approximators like Neural Networks can be used for our algorithm. Indeed, in the HJB case, the Policy Iteration algorithm leads to solving backward Kolomogorov equations, which have previously been treated with neural networks in \cite{grohs-jensen}. 


In the following section we introduce the optimal control problem and the concept of the Bellman equation.
Section \ref{sect:pol_it} is devoted to the Policy Iteration in a function space. In the subsequent sections the Least-Squares approximation and our function approximator, the Tensor Trains,  are introduced. The final section is devoted to the presentation of the numerical results, where we cover some low-dimensional problems and one problem in dimension six.

\section{Exit Time Control Problem of a Stochastic ODE}\label{sect:exit}
We assume, that the stochastic ODE in $\Xi \subset \mathbb R^n$ open, given by
\begin{align}
    d X_t &= b(X_t)dt + \sigma(X_t) d W_t + g(X_t) u_t dt \label{eq:open_loop_sde}\\
    X_0 &= x, 
\end{align}
describes the state of a system at time $t$, where 
$X_t\in \mathbb{R}^n$, $b:\mathbb{R}^n \to\mathbb{R}^n$ is the gradient of a smooth potential $U$ with bounded derivatives, $g: \Xi \to \mathbb R^{n, m}$, $\sigma:\mathbb R^n\to \mathbb R^{n, n}$ are smooth with bounded derivatives and $W$ is a Wiener process{.  Additionally,  $u_t \in \mathbb{R}^m$ is a control parameter adapted to the process $X_t$. In the context of dynamic programming $u_t$ is also known as action, decision or policy of the controller. 
For each control $u$  we can define the exit time $\eta$

$$\eta = \inf \{t > 0 | X_t \not\in \Xi \}$$

and the cost function
\begin{align*}
    {\mathcal{J}}(x,u) &= \mathbb E \Big[\int_{0}^{\eta}  c(X_t) + u_t^T B u_t dt| X_0 = x \Big],
\end{align*}
where $c:\mathbb R^n\to\mathbb R$ is some given continuous, positive function and $B \in \mathbb R^{m,m}$ is positive definite.

Our goal is it to find the control $u$ with minimal cost, i.e.
\[
\min_u {\cal J}(x,u).
\]
Here, we do not specify the space that we minimize over. Formally speaking, we consider the space of controls mapping to $\mathbb R^m$ that are measurable and adapted to a filtration induced by the Brownian motion. As we are later considering feedback controls, we omit the technical details and instead refer to \cite{Bismut1978, Soner+Fleming} and references therein.

We define the value function $v^*:\mathbb{R}^n \to \mathbb{R}$ as minimum of the cost functional over all controls, i.e.
\[ 
v^*(x):=\min_{u} \mathcal J(x, u) \text{ for all } x \in \Omega. 
\]
The control $u^*$ is optimal if $v^*(x)= \mathcal J(x, u^*)$ holds.

Under the assumption that there exists a Lipschitz continuous feedback control $\alpha:\mathbb{R}^n \to \mathbb{R}^m$ with finite cost for any initial state $x\in \Omega$
we replace \eqref{eq:open_loop_sde} by the closed loop system 
\begin{equation}\label{eq:closed_loop}
    d X_t = b(X_t) + \sigma(X_t) d W_t + g(X_t) \alpha(X_t)
\end{equation}
and denote the corresponding state by $X_t^\alpha$ to stress the dependence on the feedback.
It has been shown under suitable regularity assumptions on the right hand side of the SDE, that \eqref{eq:closed_loop} is well defined and is differentiable with respect to the intial values \cite[chapter 2]{ArnoldLudwig1998RDS},\cite{Arnold1995PerfectCT}. Since we are considering time-homogeneous Itô diffusions, our processes fulfill the (strong) Markov properity \cite{Oksendal}.

In the following we only consider feedback laws that give us finite costs. Thus, we define the 
policy evaluation function $v^\alpha$ with respect to a fixed feedback law $\alpha$ as
\begin{equation}\label{eq:v_alpha}
    v^\alpha(x) :=  \mathbb E \Big[ \int_0^\eta c(X_t^\alpha) + \alpha(X_t^\alpha)^T B \alpha(X_t^\alpha) dt | X_0 = x \Big].
\end{equation} 
\begin{remark}\label{rem:pol_evaluation}

In our algorithm we will later compute $v^\alpha$ approximately. In this case the equality \eqref{eq:v_alpha} is not given and we the distinguish between the policy evaluation function, which is $v^\alpha $ as in \eqref{eq:v_alpha}
and the policy estimation function, which is $v^\alpha$ computed by our algorithm.
\end{remark}
We define the set of feedback laws that induce finite costs 
$$F = \{ \alpha:\mathbb R^n \to \mathbb R^m\ |\ v^\alpha(x)<\infty \text{ for all }  x \in \mathbb R^n \}$$
and assume that there exists an optimal, Lipschitz continuous feedback law.
This ensures that the value function is the policy evaluation functional of the optimal feedback, i.e.
\[ v^* = v^{\alpha^*}. \]
If the value function $v^*$ is known and differentiable, the optimal control is given explicitly by \cite{DeterministicHJB}
\[ \alpha^*( x) = - \frac 1 2 B^{-1} g( x)^T \nabla v^*( x).\]

Moreover, abbreviating $ r^\alpha(x) = c(x) + \alpha(x)^T B \alpha(x)$ and $\tau \land \eta = \min \{ \tau, \eta \}$, a stochastic Bellman equation \cite{Buckdahn2016} is obeyed for every $\tau > 0$ 
\begin{align}
v^*(x) &=  \mathbb E \Big[ \int_{0}^{\eta \land \tau} r(X_t^{\alpha^*}) dt + v^*(X_{\eta \land \tau}^{\alpha^*})| X_0 = x \Big] \label{eq:Bellman1}  \\
\alpha^*( x) &= - \frac 1 2 B^{-1} g( x)^T \nabla v^*( x)  \label{eq:Bellman2}
\end{align} 
with Dirichlet boundary condition 
\[ v^*(x) = 0 \text{ on } \partial \Xi. \]
We later employ the Policy Iteration algorithm to solve this coupled equation by alternating between the value updates given by \eqref{eq:Bellman1} and the policy updates given by \eqref{eq:Bellman2}. 
In preparation to that we first notice that by fixing a policy $\alpha$ this coupled equation becomes uncoupled and a linear function equation is remaining
\begin{equation}\label{eq:lin_bellman}
    v^\alpha(x) = \mathbb E \Big[ \int_{0}^{\eta \land \tau} r(X_t^{\alpha}) dt + v^\alpha(X_{\eta \land \tau}^{\alpha})\Big]:=\mathbb E \Big[ \int_{0}^{\eta \land \tau} r(X_t^{\alpha}) dt + v^\alpha(X_{\eta \land \tau}^{\alpha})| X^{\alpha}_0 = x \Big]. 
\end{equation}
Note, that the expectation value is a conditional expectation value with respect to the initial value $X_0=x$. For the ease of notation, in the sequel we sometimes drop the condition $ X_0=x$ when the context is clear.

\section{Policy Iteration}\label{sect:pol_it}
A typical approach to solve the Bellman equation \eqref{eq:lin_bellman} would be to use a fixed point iteration in the values, with some given inital guess. This method is known as value iteration e.g. \cite[Chapter 3]{Bertsekas}.
When computing \eqref{eq:lin_bellman} with the value iteration, a sequence of functions $\{v_k\}_k$ is generated, such that $v_k \to v^*$
under suitable conditions. However, for a value iteration it must be clarified how to discretize the policy(space). 

In the following, we take a different path
by using the Policy Iteration, where a sequence of polices $\{\alpha_k\}_k$ instead of values is generated. 
It has been understood in \cite{oster2019approximating} that the policy evaluation function $v^\alpha$ in the deterministic setting can be computed via the Koopman operator \cite{Koopman315}. 
The Koopman operator is a linear transfer operator allowing to transfer a system with a non-linear evolution to a linear system in function space. The structure of the eigenvalues and eigenfunctions of the Koopman operator have been investigated to obtain a coarse grained description of the system \cite{Dellnitz_2000,DellnitzJunge,SCHUTTE1999146,Koopman}.

To do so we rewrite \eqref{eq:lin_bellman} by using the Koopman operator \citep{PerronFrobeniusKoopmanSchuette,LasotaMackey, AppliedKoopmanism, KoopmanofRandomDynSys} with a slight modification to incorporate the exit time
$$K_\tau^\alpha[v] = \mathbb{E} [ v( X_{\eta \land \tau}^\alpha)| X^{\alpha}_0 = x ].$$
This allows us, to reformulate equation \eqref{eq:lin_bellman} as operator equation
\begin{equation}\label{eq:Kop_hjb_transform1}
(I -K_\tau^\alpha)v(x) = \mathbb E[\int_0^{\eta \land \tau} r(X_t^\alpha)dt | X^{\alpha}_0 = x ].
\end{equation}
With the operator equation \eqref{eq:Kop_hjb_transform1} we can now give the Policy Iteration algorithm.

\begin{algorithm}[H]
\SetAlgoLined
\caption{Policy Iteration for solving \eqref{eq:Bellman1}}\label{algo:pol_it_continuous_global}
\SetKwInOut{Input}{input}\SetKwInOut{Output}{output}
\SetKwInOut{Output}{output}\SetKwInOut{Output}{output}
\Input{A Policy $\alpha_0 \in F$.}
\Output{An approximation of $v^*$ and $\alpha^*$.}
Set $k = 0$.

\While{not converged}{
  Solve the linear equation
\begin{equation}\label{eq:stoch_lin_hjb}
     (I-K_\tau^\alpha) v_{k}(x) = \mathbb E \big[\int_{0}^{\eta \land \tau} r(X_s^{\alpha_k}) ds | X^{\alpha}_0 = x \big],
\end{equation}
then update the policy according to
\begin{equation*}
\alpha_{k+1}(x) = - \frac 1 2 B^{-1} g(x)^T \nabla v_{k}(x) .
\end{equation*}
$k = k+1$.
}
\end{algorithm}
Note that in the algorithm we have to choose an initial policy $\alpha_0 \in F$. In some cases this is a particular hard challenge, as the policy has to lead to finite cost for every initial state $x$. However, the Wiener process ensures that the uncontrolled dynamics driven by a potential arrive at the exit set in finite time  almost surely. Thus, we initialize the Policy Iteration with the zero control $\alpha \equiv 0$.
For solving \eqref{eq:stoch_lin_hjb} in the above algorithms we face the following two problems:

    \paragraph{Linearized Bellman} In principle the point values of   the linearized Bellmann equation 
    \eqref{eq:stoch_lin_hjb}
    can be computed by Monte Carlo Methods  such as Euler-Mayurama. Let us assume that we have computed pointwise values.
    How can we 'interpolate' between these values to obtain an approximation 
    of the policy evaluation function?
    For this purpose we propose a Least Squares approach, where the nodes are sampled randomly, 
    i.e.  we use Monte Carlo integration. We call this approach {\em variational Monte Carlo}, introduced in section \ref{sect:vmc}.\
    
       \paragraph{Model Class} The equation \eqref{eq:stoch_lin_hjb} is given in an infinite dimensional space $V$ and we need a finite dimensional ansatz space, or at least a  set  $U  \subset V$ of computable aproximation $ v_{\varepsilon ,k} $ of $v_k$ to achieve a desired accuracy $\varepsilon$. We choose the intersection of a scaled ball with a submanifold $\mathcal M$. Therefore, we have to approximate the function $v_k\in U\subset \mathcal{M} \subset  V$. We propose Tensor Trains and Tree based Tensors for $\mathcal M$ in Section \ref{sect:TT} to tackle this challenge.}

\todo{mathias aenern und alle dann lesen}
\section{Variational Monte-Carlo}\label{sect:vmc}

We now elaborate on how to tackle  the computational bottleneck, i.e. 
the linearized Bellman equation \eqref{eq:stoch_lin_hjb},
\begin{align}
     (I-K_\tau^\alpha) v_{k}(x) &= \mathbb E \big[\int_{0}^{\eta \land \tau} r(X_s^{\alpha_k}) ds | X^{\alpha}_0 = x \big]. \\
      \iff v_k(x) &= \mathbb E \big[\int_{0}^{\eta \land \tau} r (X_s^{\alpha_k } ) ds + v_k (
X_{\eta \land \tau}^\alpha
) | X^{\alpha}_0 = x  \big]. \label{eq:stoch_lin_fixedpoint}
\end{align}
This is done in three major steps. First, we interpret the above equation as a fixed-point equation, then we formulate the subproblems as a Least Squares problem on a finite dimensional function space and, finally, we use Monte-Carlo quadrature to integrate within the state and probability space.

The policy $ \alpha:= \alpha_k$ is given and we  assume that $\hat v$  on the r.h.s. of the equation below is given as well.
%
Let us introduce $\tilde v \in V $ such that 
\begin{equation}
\tilde v(x) := \mathbb E \big[\int_{0}^{\eta \land \tau} r (X_s^{\alpha } ) ds + \hat v (
X_{\eta \land \tau}^\alpha
) | X^{\alpha}_0 = x  \big] \ , \ x \in \Omega ,     
\end{equation} 
Note that if $\tilde v = \hat v$ we have found a solution to \eqref{eq:stoch_lin_fixedpoint}. As this equation is posed in an infinite dimensional function space, we first rewrite it as a Least-Squares problem on a finite dimensional subspace.

More exactly, consider the Hilbert space $ V := L^2 (\Omega ) $, (or more generally $ V := L_2 ( \Omega , \rho ) $ with some probability density $ \rho$) such that $\hat v, \tilde v \in V$. Then, we have
\begin{equation}\label{eq:Riskfunctional1}
\tilde v = \argmin_{v \in V} \mathcal R^\alpha(v, \hat{v} ), \quad \mathcal R^\alpha(v, \hat{v} ) = \| v(\cdot )  -\mathbb E \big[\int_{0}^{\eta \land \tau} r(X_s^{\alpha }) ds + \hat v(
X_{\eta \land \tau}^\alpha
) | X_0^\alpha  = \cdot \big] \|_V^2   
\end{equation} 
and since  $\tilde v \in V$, we have $\mathcal R^{\alpha}(\tilde v, \hat v ) = 0$. 
Indeed,  we are seeking a solution  $ \hat{v} = \tilde{v}$ which constitutes a fixed point problem 
\begin{equation}\label{eq:Riskfunctional2}
\tilde v = \argmin_{v \in V} \mathcal R^\alpha(v,  \tilde{v}), \quad \mathcal R^\alpha(v,  \tilde{v}) = \| v -\mathbb E \big[\int_{0}^{\eta \land \tau} r(X_s^{\alpha }) ds + \tilde v(
X_{\eta \land \tau}^\alpha
)| X_0^\alpha  = \cdot  \big] \|_V^2.
\end{equation}
However, finding the exact  minimizer $\tilde v \in V$  is infeasible,  and thus we further restrict to a finitely representable  compact subset 
$U \subset V$.

Classically, $U \subset V$ is a   closed ball of some finite dimensional subspace $V_p\subset V$. However, in many applications the subspace $V_p$ is high-dimensional, which makes computations impracticable. Thus, we introduce a lower dimensional submanifold $\mathcal M\subset V_p$ and consider  $U \subset \mathcal{M} \subset V_p$ to be a compact subset
feasible for computational treatment,  having an intrinsic data complexity, which can be handled
by our technical equipment. In our case $\mathcal{M}$ is the set of tensor trains of bounded (multi-linear) rank
$\mathbf{r} = (r_1, \dots, r_n)$ and $U$ is the set of rank $\mathbf{r}$ tensors with uniformly bounded norm embedded in the space $\mathcal{M}\subset V_p$ of tensor product polynomials of multi-degree ${\bf p}$, which will be covered in detail in the following section.
For $ U \subset \mathcal{M} $ the  Least-squares approximation is defined by the minimizer \eqref{eq:Riskfunctional2}
\begin{equation}\label{eq:minimizer-risc}
    v^{\alpha}_U = \argmin_{ v \in U} \mathcal{R}^{\alpha}(v,v ).
\end{equation}
However, the numerical treatment of the above 
minimization problem  \eqref{eq:minimizer-risc} is still infeasible, due to the presence of the high-dimensional integrals over $\Omega \subset \mathbb{R}^d$. To handle this problem  we replace the exact integral by a numerical quadrature. 

We compute the norm $ \| . \|_V$ using Monte-Carlo integration, e.g.
\begin{align}
\begin{split} \label{eq:emp-risk}
    v^{\alpha}_{(U, N, M) } &= \argmin_{ v \in U} \mathcal{R}_{N, M}^{\alpha}(v,v) , \\
  \mathcal R_{N, M}^\alpha(v,v) &=  
  \frac{1}{N} \sum_{i=1}^N \big|v(x_i)- \frac 1 M \sum_{j = 1}^M\int_{0}^{\eta \land \tau} r (x_i^{j,\alpha} (t)  ) dt + v
  ( x_i^{j,\alpha } ({\eta \land \tau})
  )\big|^2,
  \end{split}
\end{align}
where $x_i^{j, \alpha}  (0)  = x_i$ for $i = 1, \dots, N$. 
We remark, that we have two Monte Carlo approximations: The first for integrating the stochastic differential equation \eqref{eq:closed_loop}
with different paths $ t \mapsto x^{j, \alpha } (t)$, $ j=1 , \ldots , M$.
And the second Monte Carlo integration for different
 initial values $x=x_i$, $i=1, \ldots N$ for setting up the Least Squares functional \eqref{eq:minimizer-risc}. 
The integral term $\int_{0}^{\eta \land \tau} r (x_i^{j,\alpha} (t)  ) dt $ in \eqref{eq:emp-risk} is then computed by a trapezoid rule.
The different paths are computed by Euler Mayurana scheme \cite{KloedenPlaten} at discrete times $ t_k$ and the remaining integral in the formulars above are approximated 
by trapezoidal rule.

Let us highlight that the  the 
input data  are noisy due to stochastic nature of the SDE. 
Therefore, the Least Squares method is prone to over-fitting problems. 
Moreover, an accurate computation of the updated policy, and therefore an accurate approximation of the gradient of $v^{\alpha } $  is 
ultimately important for the convergence of the Policy Iteration. 

\paragraph{Regularisation}

To improve the accuracy, we can enforce better regularity of $v$ by choosing an appropriate norm $\| \cdot \|_F$.
In the numerical calculations,  we add a regularization term, such that the actual 
risk functional is
\begin{equation}
    \label{eq:loss_funct}
    \tilde {\mathcal R}_{N, M}^\alpha (v,v) : = \mathcal R_{N, M}^\alpha(v,v) + \delta \| v \|_F^2.
\end{equation}
By choosing $\| \cdot \|_F$ as the $\ell_2$ norm of the coefficient tensor $ A $ of $v$ the regularization term depends on the choice of univariate basis functions.
Presently, we used mixed (tensor product) Sobolev norms and refer to remark \ref{rem:ansatz_fun} for a brief discussion.
To avoid deviation of the solution caused by the penalty term,
we decrease the penalty parameter adaptively
during the ALS iteration process in dependence of the current residual \cite{ALS}.
For any other iterative solver of \eqref{eq:emp-risk} this procedure can be done analogously.

We have experienced that this part plays a crucial role for the performance of the algorithm.  In our present examples,   our method provided  quite accurate results. However, for non-smooth viscosity solutions arising from more difficult problems  we expect 
that improved techniques will be required. 

The present optimization  problem is tractable by local optimization methods on non-linear manifolds. 
Nevertheless, it remains hard to find an exact minimizer,
see e.g. \cite{Bachmayr-Uschmajew-Schneider} for further discussion.

\paragraph{Error Estimates} The theoretical  justification 
of this Least Squares Monte Carlo approach  
is in a very early stage. 
Indeed,  we are committing variational crimes,  since we have replaced the original risk functional by an empirical risk functional. 
This introduces an additional error term,  even if we assume that we 
have computed the exact minimizer of \eqref{eq:HJB-opt}.  
For first theoretical results, we 
import well known results from {\em empirical risk minimization} in  machine learning \cite{cucker-smale,steinwart}. 
Empirical risk minimisation has been considered 
for the regression problem in statistical learning theory. 
However, the present problem is not directly a regression problem, 
but the theory  \cite{cucker-smale} can be straightforwardly extended
to the present optimization problem, which we called 
{\em Variational Monte Carlo}.  This term  has been invented in physics earlier, but Monte Carlo Least Squares method seems be 
also an appropriate name. 
The error 
$\mathcal{E} = \| v^{\alpha} - v^{\alpha}_{(U,N,M)} \|_v^2 $
is split into three parts 
$$\mathcal{E} = \mathcal{E}_{\text{approximation}} + \mathcal{E}_{\text{generalization}} + \mathcal{E}_{\text{optimization}}   . $$
Due to the uncertain nature of the problem, we cannot expect to show convergence for the generalization term. Instead, we consider convergence in probability. In particular, the probability that a given error estimates 
fails, i.e. $ \mathbb{P} [ \mathcal{E}_\text{generalisation} > \varepsilon ] $,  decays exponentially with the number of sample points.

As a first result we recall to following corollary.
Under certain assumption, one can show that



\[ \mathbb P[\|v^\alpha-v^\alpha_{( U, M,N)}\|_{L^2(\Omega)}^2> \varepsilon] \leq c_1(\varepsilon, U) e^{-c_2 N \varepsilon^2} \]
with $c_1(\varepsilon, U), c_2>0$.

See e.g Theorem 4.12, Corollary 4.19 and Corollary 4.22 from \citep{VMC}. 

Overfitting  effects introduced by the above interpolation
procedure 
spoil the computation of the optimal policy  $ \alpha$ more dramatically, since this requires the gradient of $v^{\alpha}$ \eqref{eq:HJB-opt}. We have experienced this effect in our computation. 
For a theoretical justification,
it would be desirable that the  error 
can be estimated in much stronger norms, e.g   w.r.t.   Lipschitz-norms etc..
In the result mentioned above,  it was  only measured w.r.t an $L_2$-norm.

\paragraph{A Regression Problem}
Let us highlight, that the present approach is NOT {\em learning}.   It is NOT a statistical task,  but a numerical method to solve an operator equation. 
Instead of solving \eqref{eq:emp-risk} directly, we use a fixed point iteration
$$ 
v_{k+1} :=   \argmin_{ v \in U} \mathcal{R}_{N,M} ^{\alpha}(v,v_k ).
$$
In fact, this constitutes a regression problem.

\section{Tree Based Tensor Representation - Tensor Trains}\label{sect:TT}

For large  dimensions $n$,  traditional ansatz functions, e.g. finite elements, splines, multi-variate polynomials etc. are not appropriate  for the numerical solution of the PDE, since they  are facing the curse of dimensions. 




To this end we choose an underlying finite dimensional but large subspace
$ {V}_p \subset {V}= L^2(\Omega)$ 
for the approximation of the sought value function. 

First we 
choose a suitable  approximation space for univariate approximation of functions $x_i \mapsto f(x_i)$, $i=1, \ldots, n$. 
Presently, we have taken
 one-dimensional polynomials $\phi_i = \text{span} \{ \phi_{i_l} | 0 \leq i_l \leq p_i \}$ of degree $p_i$. 
 However, other choices like splines waveletes etc. are also possible. 
 
 For the $n$-variate case, we consider the tensor product of such polynomial spaces
$$ {V}_p :=
\mbox{span }  \{ \phi_1 \otimes \cdots \otimes \phi_n : \mbox{deg} \phi_i \leq p_i 
\} .
$$
This is a space of multivariate (tensor product) polynomials with bounded multi-degree $ \mathbf{p} = ( p_1, \ldots, p_n)$.
For the sake of simplicity we  have chosen  the same degree in all coordinates, i.e.$ p_i = p$, $i=1, \ldots , n$.

A function $ q \in V_p$ can be expanded w.r.t. to tensor product basis functions via
$$ 
q ( x_1 , \ldots , x_n ) = 
\sum_{i_1 , \ldots , i_n = 1}^{p_1, \dots, p_n}
A ( i_1, \ldots , i_n ) \phi_{i_1} (x_1) \cdots 
\phi_{i_n} (x_n) .
$$
Interpreting the coefficient representation  
$ (i_1, \ldots , i_n) \mapsto A ( i_1, \ldots , i_n )  $
of a polynomial $ q$ in this vector space as a tensor of order $n$, 
we need storage in $\mathcal O (p^n)$ for the coefficient tensor $ A \in \otimes_{j=1}^n \mathbb{R}^{p} $.

Let us note that $ \bigcup_{p\in \mathbb{N}} V_p $
 is dense in $V$. 
 Although the dimension of $ V_p$ is finite 
$$
\mbox{dim} {V}_p = (p+1)^n 
$$
it is prohibitively
large.

In the  ambient  space ${V}_p$, we consider a non-linear, possibly low-dimensional  
manifold,  given by 
{\em tree based tensor representations} (hierarchical (Tucker) tensors - HT tensors) \cite{Hackbusch-buch}.
In the present applications,  we choose 
so-called 
tensor trains (TT tensors), invented by Oseledets in \citep{Oseledets,Oseledets2},
which has {considerably} smaller dimensions \cite{TTTensor}. 
They have been applied to various high-dimensional PDE's \citep{Khoromskij-book}, but the parametrization has been used in quantum physics much earlier 
as {\em Matrix Product States} and {\em Tensor Network States}, successfully for the approximation of spin systems and Hubbard models. 
For a good  survey  we refer to the papers 
\citep{Hackbusch2014,Bachmayr-Uschmajew-Schneider,Legeza-Schneider,Hackbusch-Acta}. 
The tensor train representation  have appealing properties
making them attractive for treatment of the 
present problems, compare \citep{TensorKunisch}. For example they contain sparse polynomials, 
but are much more flexible at a price of a slightly larger overhead, see e.g. \citep{Dahmen-Bachmayr} for a comparison concerning parametric PDEs. 
Let us  give brief introduction for a first understanding.

 In order to get some notion of the representation and  compression, we introduce the TT-rank $\mathbf{r} \in \mathbb N^{n-1}$ of the tensor $A\in \mathbb R^{(p_1,\dots,p_{n-1})}$ as element-wise smallest tuple such that
\[ A(i_1,\dots,i_n)= \sum_{k_1,\dots,k_{n-1} = 1}^{r_1,\dots,r_{n-1}} U_1(i_1,k_1)\cdot U_2(k_1,i_2,k_2)\cdots U_n(k_{n-1},i_n) \]
holds for some $U_i\in \mathbb R ^{r_{i-1},p_i,r_i}$ for $i=1,\dots,n$.
The TT-rank is well defined and the set of  tensors of fixed TT-rank $\mathbf{r}$ 
forms a smooth manifold of dimension in $\mathcal O (npr^2)$ \citep{TTTensor} in contrast  to  $ \mathcal{O} (n^p)$ of the ambient linear space ${V}_p$. 
Taking the closure of this set,  see e.g. \cite{Hackbusch-buch} 
we allow also tensors with smaller TT-rank denoted by $\mathcal M := \mathcal{M}_{\leq \mathbf{r}} $ \cite{Bachmayr-Uschmajew-Schneider}. This slightly larger set forms 
an 
{\em algebraic variety}
\citep{KUTSCHAN,Bachmayr-Uschmajew-Schneider}. 
However,  numerical routines  like ALS \cite{ALS} 
do not differentiate between the variety and the manifold.
For a survey and mathematical theory we refer to the literature, e.g. 
\cite{Hackbusch-buch,Hackbusch-Acta,Bachmayr-Uschmajew-Schneider}.


TT tensors can represent polynomials as follows. Let us consider the vectors
\begin{align*}
P_i\colon \mathbb R \to \mathbb R ^{p_i+1} \text{ with } P_i(x) = \begin{bmatrix}1 \\ \phi_0 (x)  \\ \phi_1 ( x) \\ \vdots \\ \phi_{p_i } (x ) \end{bmatrix} =\begin{bmatrix}1 \\ x \\ x^2 \\ \vdots \\x^{p_i}\end{bmatrix}.   
\end{align*}

Then

\begin{align*}
p(x_1,\dots,x_n)  & = \sum_{i_1,\dots,i_n}^{p_1,\dots,p_n} \sum_{k_1,\dots,k_{n-1}}^{r_1,\dots,r_{n-1}} U_1(i_1,k_1) U_2(k_1, i_2, k_2)\cdots 
\cdot U_n(k_{n-1},i_n) \big(P_1(x_1)\big)_{i_1}(P_2(x_2)\big)_{i_2}\cdots \big(P_n(x_n)\big)_{i_n} \\
& = \sum_{i_1,\dots,i_n}^{p_1,\dots,p_n} \sum_{k_1,\dots,k_{n-1}}^{r_1,\dots,r_{n-1}} U_1(i_1,k_1) U_2(k_1, i_2, k_2)\cdots 
\cdot U_n(k_{n-1},i_n) \phi_{i_1} (x_1) \phi_{i_2} (x_2)
\cdots \phi_{i_n} (x_n)  
\end{align*} 
is a multivariate polynomial
of degree $\left(\sum_i p_i\right)$. 

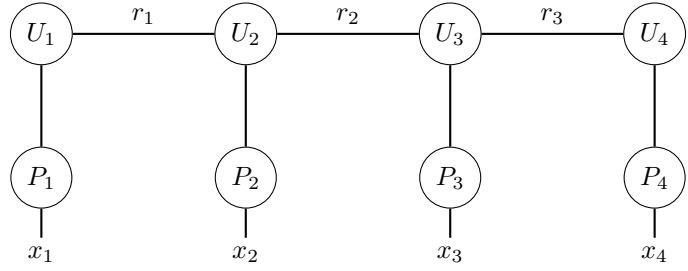
\begin{figure}
    \centering
    \begin{tikzpicture}

\begin{scope}[every node/.style={scale=1,minimum size=5mm,circle, draw,  fill=white}]
    
    \node (A1) at (0,0) {$U_1$};

   	\node [ position=0:2.3 from A1](A2) {$U_2$};  	
   	\node [ position=0:2.3 from A2](A3) {$U_3$};   	
   	\node [ position=0:2.3 from A3](A4) {$U_4$};

	\node [ position=-90:1.5 from A1](B1) {$P_1$};
   	\node [ position=-90:1.5 from A2](B2) {$P_2$};
   	\node [ position=-90:1.5 from A3](B3) {$P_3$};
   	\node [ position=-90:1.5 from A4](B4) {$P_4$};

\end{scope}

	\node (s) at (0,1) {};

	\node [ position=-90:.8 from B1](C1) {$x_1$};
   	\node [ position=-90:.8 from B2](C2) {$x_2$};
   	\node [ position=-90:.8 from B3](C3) {$x_3$};
   	\node [ position=-90:.8 from B4](C4) {$x_4$};

\begin{scope}[every node/.style={}]

\end{scope}

\begin{scope}[every edge/.style={draw=black,thick}]
	\path [-] (A1) edge node[midway,left] [above] {$r_1$} (A2);
	\path [-] (A3) edge node[midway,left] [above] {$r_2$} (A2);
	\path [-] (A3) edge node[midway,left] [above] {$r_3$} (A4);

	\path [-] (A1) edge node[midway,left]  {} (B1);
	\path [-] (A2) edge node[midway,left]  {} (B2);  
	\path [-] (A3) edge node[midway,left]  {} (B3);  
	\path [-] (A4) edge node[midway,left]  {} (B4);

	\path [-] (B1) edge node[midway,left]  {} (C1);
	\path [-] (B2) edge node[midway,left]  {} (C2);  
	\path [-] (B3) edge node[midway,left]  {} (C3);  
	\path [-] (B4) edge node[midway,left]  {} (C4);  

\end{scope}

\end{tikzpicture}
    \caption{Graphical representation of TT tensor train induced polynomial in four variables.}
    \label{TT_polynomial}
\end{figure}

Using the graphical tensor network representation \cite{Bachmayr-Uschmajew-Schneider,Legeza-Schneider} this polynomial can be interpreted as in Figure \ref{TT_polynomial}.

\begin{remark}\label{rem:ansatz_fun}
Note that other  polynomial basis functions  can be chosen as well. For numerical reasons we choose a set of orthogonal polynomials, e.g. $ \phi_i = \ell_i$ Legendre polynomials. 
In this case, Parseval formula provides  a norm equivalence between  $ L_2 $ and $\ell_2$, which guarantees stability of our representations
and approximation schemes. Presently,  we have chosen 
one-dimensional $H^1 (I)$ orthogonal polyomials. The stability is enforced by an additional regularization term, and the penalty parameter has been adaptivley reduced through the iteration procedures. 
This procedure enforces the approximations to have small 
$L_{\infty}$ and even Lipschitz bounds. 

In general,  the set of one-dimensional  basis functions can be modified to fit better to other norms rather than $L_2$ or $H_1$.
\end{remark}

It turns out, that optimization procedures in this TT format can be solved by consecutively optimizing one component $U_i$  while the others are fixed. This alternating Least-Squares (ALS) algorithm converges to a local minimum \citep{ALS}. Further details on the implementation in a similar context can be found in \cite{oster2019approximating}.

\begin{remark}
The present tensor ansatz has been proved by our experience to provide an advantageous choice, however there are some  alternatives, well known in machine learning 
which can be used for the present purpose in same fashion or  with some more or less obvious modifications.
Among them are 
sparse grids \cite{bungartz_griebel_2004}, sparse polynomials   \cite{BachmayrCohenSchwab17} , 
kernel methods (SVM) \cite{smola}, in  particular with polynomial kernels,
and deep neural networks (DNN) \cite{Goodfellow-et-al-2016} . 

In this respect,  the essence of the present paper is not restricted to tree based tensor methods. 
\end{remark}

\section{Formal Scaling with Respect to the Spatial Dimensions}
We 
add a brief discussion about the computational  complexity, and how the computational complexity scale with the spatial dimensions $n$, and how the HJB is prone to the curse of dimension.

Let us assume that we want to achieve a fixed accuracy $\epsilon$, i.e. we do not consider the scenario $\epsilon \to 0$. This is motivated because
we want to keep the feedback law fairly simple, since this is required for an online feedback law. 

The number of degrees of freedom of the underlying TT tensor is $ K \leq  n(p+1)r^2 = \mathcal{O} (n)$  for fixed accuracy. 
Note that $p=p(\epsilon) $, and $r= r(\epsilon) $ and will kept as constants in the sequel. 
(Presently $ p=10,n=6,r= 5 $.)  
In this regime $ K \sim n $ scales linearly with $n$ instead of exponentially. This linear scaling 
behavior for storing the value function,  seems to be quite optimal. We have rendered the curse of dimensions in a perfect way.

We further assume that we need at least 
$N= \mathcal{O} (K) $ sample points, which is  very optimistic. This is the best scaling we can expect, and extremely optimistic, and cannot be improved by other methods like kernel methods or DNN. 
The best  proven rate for linear Least Squares methods is $ N \sim K \log K$ \cite{cohen:hal-01354003}, and we neglect further logarithmic terms.  
Therefore the total numerical Work scales at least $ \mathcal{O} (n N)  = \mathcal{O} (n^2 ) $.
With the present approach we have to calculate at each sample points $ M $ paths. 
Then, the 
total numerical work is $\mathcal{O} (n^2 M S)$, where $S$ is the work for computing a single path.

For linear  function $f$ fully connected, i.e. it is represented  by fully populated matrix,
the minimal cost for each path is 
$S = \mathcal{O} ( n^2 T )$, where $T$ is the number of time steps.
In this case we have assumed that the evaluation of the feedback law which scales with $ \mathcal{ 0 } =  (n^2 m) $, where $m$ is the number of controls. 
When $n \sim m$ the we can have an additional factor $ n$, in the scaling of $S$. This does not happen in present case. 
We summarize that the total work is 
$$ \mathcal{O} (n^4 p^2r^4 M T)  \mbox{  or for fixed accuracy } \ 
 \mathcal{O} (n^4 ).  $$ 
 In the deterministic case we save the factor $M$, since we need only a single trajectory for each sample $x^i$.i.e. for each initial condition.  
 Let us remark,
if we use the (linearized)  HJB directly, in  a Least Squares setting, 
the factor $MT$ is no longer apparent. We save also an additional factor $n$. Here, the scaling will be 
$ \mathcal{O} (n^3 p^2r^4 ) $
Indeed, this  can save computing time at a price of less stability and a loss of accuracy.  We will discuss this issue in the outlook.

In the subsequent numerical tests, we had around K= 800 DOF in our models set. We took
$N = 10 K \approx 8.000 $ sampled initial values. For each initial value $x_i$
we consider 100 paths, i.e. in each iteration step we performed $ 8 \cdot 10^5$ runs of the Euler Mayurama scheme with $100$ time steps. This part was by far the most time consuming part.
However it can be perfectly parallelized, which has not been done so far.

The above  scaling is estimated in a very optimistic fashion, and can be considered more as a lower bound. However,
there are situations where the scaling is better, e.g. the matrix representation of $f$ is sparse. It may be that multi-level Monte Carlo 
can provide an additional  better scaling. All possibilities to  reduce the present scaling  $\mathcal{O} (n^4 ) $ 
have to  be considered in next future.

For other 
model sets, known in machine learning, the scaling can be worse.
For a fully connected DNN with $\mathcal{O} ( n)$ neurons in each layer, and fixed depth $L$ we have 
$ K= \mathcal{O} (Ln^2 )  $ DOFs, and assuming $ N\sim K $ like in the above setting, 
we obtain the scaling 
$ \mathcal{O} ( n^4 L^2 n^2  M T )  = \mathcal{O} (n^6)$!
For kernel methods $ \mathcal{O} (N^3 n^2 M T ) $, where the
scaling w.r.t to the dimension $n$ number of samples $N=N_K$ 
is not clear. Assuming $ N \sim N$ seems to be quite optimistic.

\section{Numerical Results}
We present results of numerical tests for different optimal control problems. For the implementation of the tensor networks we use the open source \code{c++} library \code{xerus} \citep{xerus}. 
We also make use the python packages \code{Numpy}\citep{van2011numpy,oliphant2006guide} and \code{scipy}\cite{SciPy-NMeth}.
The calculations were performed on a AMD Phenom II 4x 3.20GHz, 16 GB RAM Fedora 31 Linux distribution. In every test we consider a compact set $E = \Xi^C$ where we want to steer the state to and a cost functional of the form
\begin{equation}\label{eq:cost}
    \mathcal{J}(x_0,u) = \mathbb E \int_{0}^{\eta} 1 +  \frac 1 2 |u_t|^2 dt.
\end{equation}
The equations are defined on a set $\Omega $ and we denote by $E = \Omega \setminus \Xi$ the set that we aim to steer the state to.
The first two tests are simple one-dimensional problems, where the exact solution is known either analytically or numerically. The third test has a two-dimensional state space and finally we test the algorithm on a $6$ dimensional state space.
\begin{remark}
In the following, we distinguish between the policy $\alpha$, the corresponding policy estimation function $v$ and the policy evaluation function $\mathcal J(\cdot, \alpha(\cdot))$. For fixed $x$, we obtain $v(x)$ by simply evaluating $v$. Here, no trajectory has to be computed. We obtain $\mathcal J(x, \alpha(x))$ by numerically integrating along the trajectory with initial condition $x$. Note that $\mathcal J(x, \alpha(x))$ is basically the numerical approximation of the cost functional with respect to a policy, defined in \eqref{eq:v_alpha}.
\end{remark}
\begin{remark}\label{rem:constants}
Within the test cases we specify the constants that we chose. Namely, the length of the trajectory $\tau$, the number of spacial samples $N$, the number of repetitions for every sample $M$ and, as we set $N$ proportional to the degrees of freedom of the tensor train representation of $v$, we also state the number of degrees of freedom.
Note that for the numerical tests the length of the trajectory does not necessarily have to be the step-size of the Euler-Majurama scheme for solving the SDE. In fact, for every test we use a step size of $0.001$ for the Euler-Majurama scheme. The length of the trajectory is mostly set to $0.1$, which means that $100$ steps within the SDE solver are used, c.f. \cite[Section 6]{oster2019approximating}.

For every test we set the regularization constant $\delta$ to be adaptive. In the beginning of every 'left-to-right sweep' within the ALS algorithm, we set $\delta$ to be the current residuum $\mathcal R^{\alpha}_{N, M}$.
\end{remark}
\subsection{Test 1. One-dimensional exit time problem: Eikonal equation}
We first test our algorithm with a simple one dimensional, deterministic problem, where the exact solution is known, namely the Eikonal equation on $\Omega = [-2, 2]$, i.e.
\[ \dot x = u, \quad x(0) = x_0, \quad \Xi = [-2, 1), E = [1, 2] . \]
Note that this problem fits into our setting by setting $\sigma = 0$ and $b = 0$. Here, the value function has the form
$$v^*(x) = \sqrt{2}(1-x),\quad x < 1.$$
\begin{figure}[h]
    \centering
    \includegraphics[width = \linewidth]{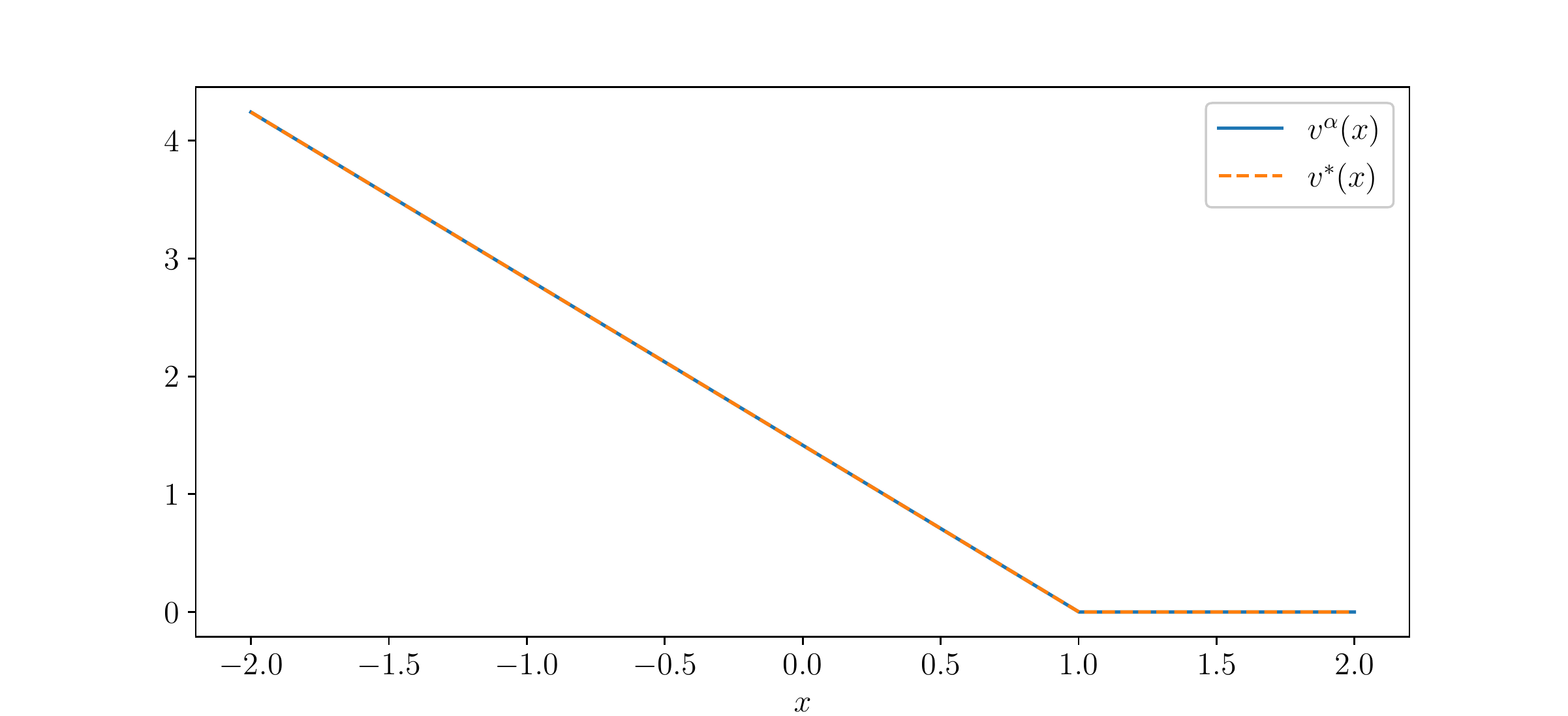}
    \caption{The value function $v^*$ and the numerical approximation $v^\alpha$.}\label{fig:1_dim_v_det}
\end{figure}
Indeed, by choosing a polynomial of degree $1$ as ansatz space, the number of samples $N = 2$ and the number of repeated samples $M = 1$, we are able to recover the value function nearly exactly, as seen in Figure \ref{fig:1_dim_v_det}. The length of the trajectory is set to $\tau=0.001$
Note that this example is particularly easy to calculate, because the optimal value function $v^*$ on the domain $(-\infty , 1) $ 
is already contained in the ansatz space, and we have set $ v(x) =0$ on the exit set $E$. 
However, it is possible to extend the domain of the ODE to $(-\infty, \infty)$ while maintaining $E = [1, 2]$. In this case,
the corresponding HJB is the well known eikonal equation 
\cite{DeterministicHJB}
$$ | v^{'} (x) |^2 =2  $$  
with boundary conditions $ v (1) = v(2) =0 $. This Dirichlet problem has multiple weak solution, the  mentioned value function $ v(x) = \max \{ 0 , \sqrt 2 (|1.5 - x| -0.5) \} $ 
is the unique viscosity solution \cite{crandall1983viscosity, DeterministicHJB}. In this respect, the present example is not so trivial and can be found in the literature for motivating the 
notion of viscosity solution. Let us remark,
that  the present  value  function $v $ is still analytic in the exterior of the target set $ E  = [1,2]$. As long as the abstract Policy Iteration (Algorithm \ref{algo:pol_it_continuous_global}) converges to a 
viscosity solution, we are going to approximate a viscosity solution. We refer to  \cite{Jacka,Kerimkulov} for a detailed discussion of this issue. We do not elaborate on this difficult issue and consider mainly  classical solutions. 
Let us remark, that in the present example $\sigma =0$ and the dynamical system is deterministic. 

We next analyze a more involved example, while still being in $1$ dimension, namely the classical double well potential.

\subsection{Test 2. One-dimensional double well potential}
\begin{wrapfigure}{r}{8cm}
    \centering
       \includegraphics[width=\linewidth]{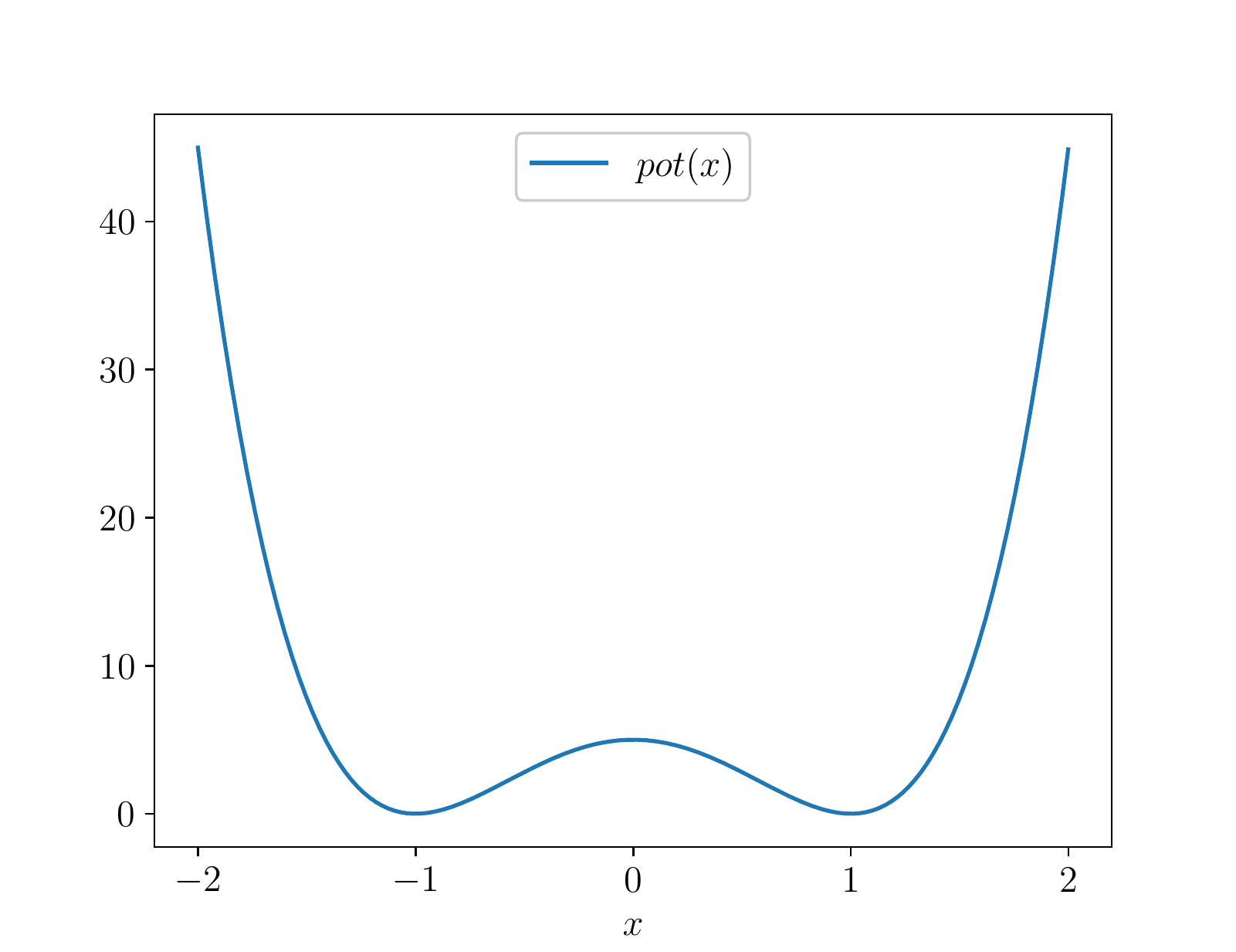}
    \caption{The double well potential $pot$.}\label{fig:1d_pot}
\end{wrapfigure}
\begin{figure}[ht]
     \subfloat[The approximated value functions and the reference solution. \label{fig:1d_values}]{%
       \includegraphics[width=0.49\textwidth]{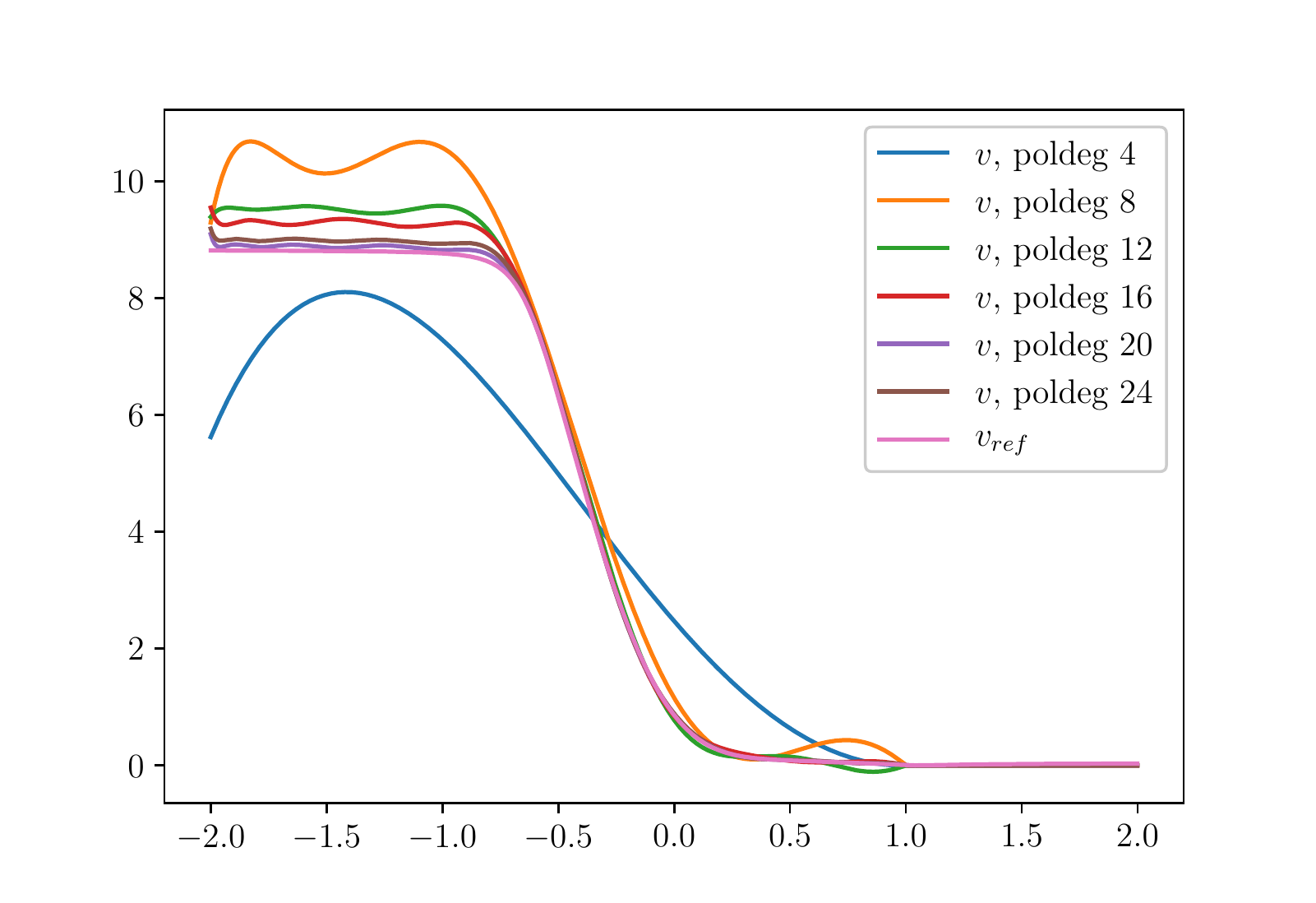}
     }
     \hfill
     \subfloat[Estimated cost and average cost by averaging over $10000$ trajectories with initial value $x_0 = -1$. Avg. cost for polynomial degree $4$ is $34.03$. \label{fig:1d_bars}]{%
       \includegraphics[width=0.49\textwidth]{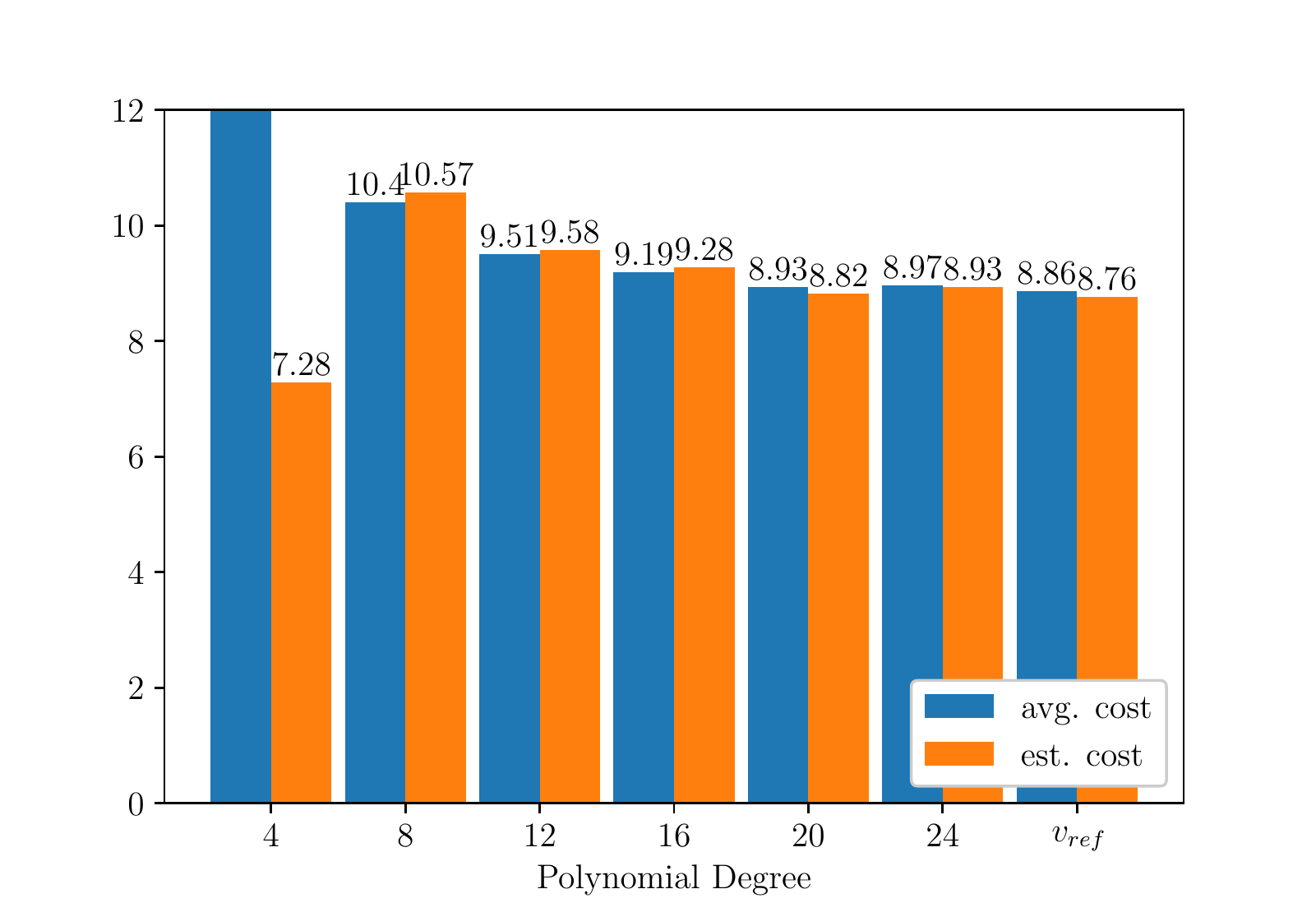}
       }
       \caption{One dimensional double well potential}
   \end{figure}

We next consider the double well potential on $\Omega = [-2, 2]$ with $\Xi = [-2, 1)$, E = [1, 2].
\[ pot(x) = 5 (x^2-1)^2, \]
visualized in figure \ref{fig:1d_pot}. 
The corresponding SDE is
\begin{equation*}
    d X_t = \nabla pot (X_t)dt + \sigma(X_t) d W_t + g(X_t) U_t(X_t)dt,
\end{equation*}
Here, we cannot expect the value function to be included in our ansatz space. Thus, we experiment with different polynomial degrees, visualized in figures \ref{fig:1d_values} and \ref{fig:1d_bars}. For the computation of the controllers we set the number of samples to $N = 10 \cdot dof$, where $dof$ is the polynomial degree increased by 1. We set the length of the trajectory to $\tau = 0.1$ and compute for every sample $M = 1000$ trajectories.
Note that the length of the trajectory consists of $100$ individual steps in the Euler-Majurana scheme, i.e. $0.1 = 100 \cdot 0.001$.

We compare the results to a reference solution that is obtained by solving the HJB equation with a finite differences scheme with $3000$ grid points.

We observe, that higher polynomial degrees yield a better approximation of the reference solution, with polynomial degree of $20$ yielding the best results. We also deduce from figure \ref{fig:1d_values} that we do not exactly reproduce $v_{ref}$.
From figure \ref{fig:1d_bars} we deduce that the performance of the controller of polynomial $20$ is less than $1\%$ higher than the performance of the reference solution.
\subsection{Test 3. Two dimensional three-hole potential}
\begin{figure}[htbp]
     \subfloat[The potential $pot$ and the exit set $\Xi^C$]{%
       \includegraphics[width=0.45\textwidth]{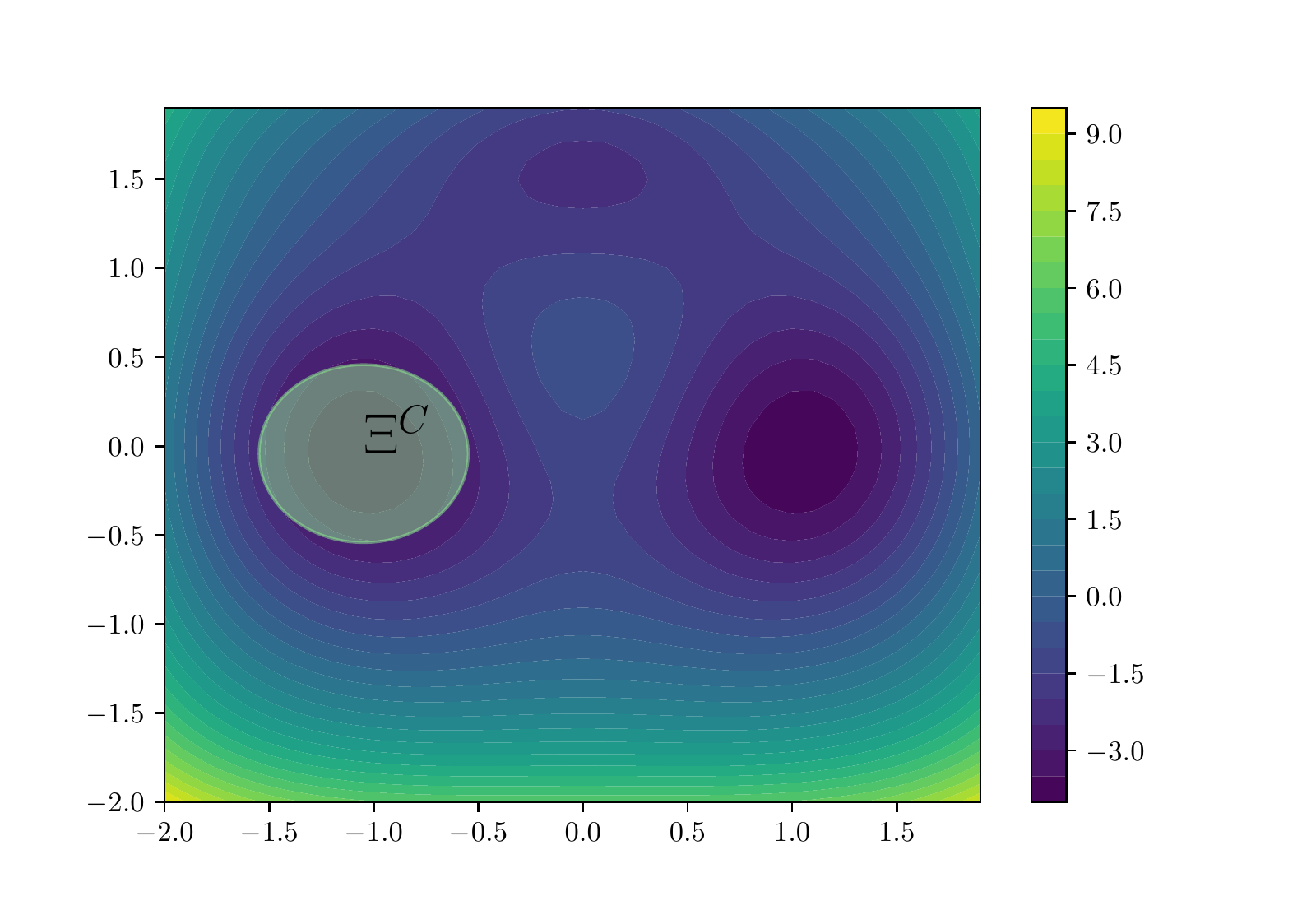}\label{fig:pot_E} 
     }
     \hfill
     \subfloat[Estimated cost and average cost by averaging over $10000$ trajectories with initial value {$x_0 = [1.8, 1.8]$}. \label{fig:2d_bars}]{%
       \includegraphics[width=0.45\textwidth]{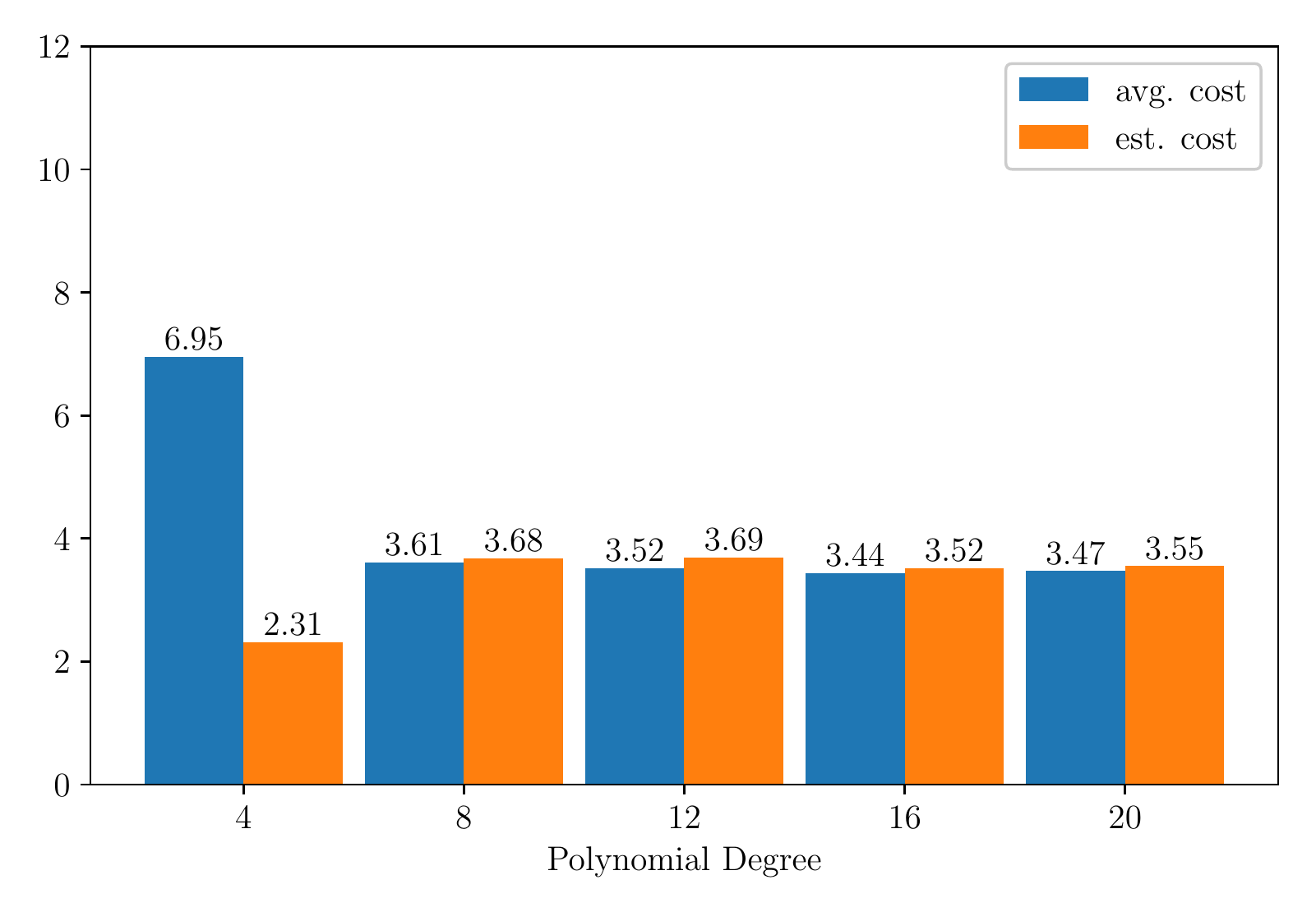}
     }
    \caption{Left: The potential $pot$. Right: Performance of controllers of different polynomial degree. }
   \end{figure}

\begin{figure}[htbp]
    \subfloat[Visualization of $10$ trajectories, uncontrolled]{
       \includegraphics[width=0.45\textwidth]{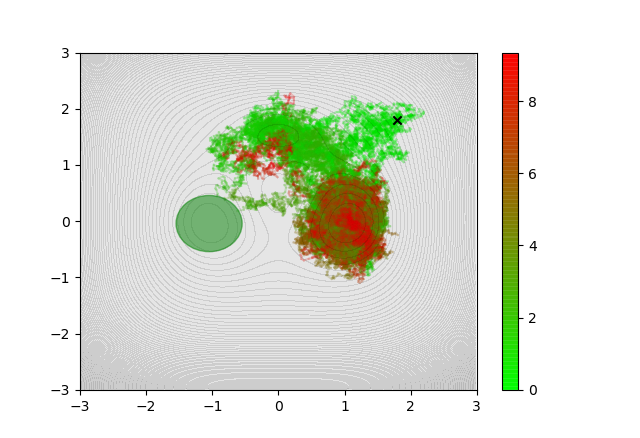}
     }
     \hfill
     \subfloat[Visualization of $10$ trajectories, controlled]{%
       \includegraphics[width=0.45\textwidth]{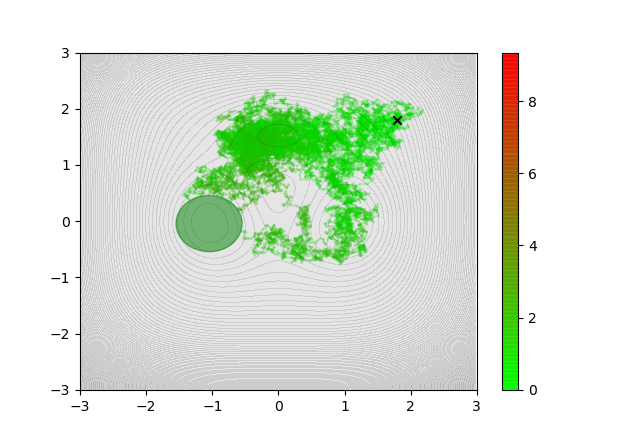}
     }
    \caption{Top: Plots of $10$ trajectories, starting at $x_0  = [1.8, 1.8]$. The colors indicate the time within the trajectory, starting with green dots and ending with red dots. We stop the simulation at $T = 10$ or when the exit set is reached. As no red dots appear for the controlled dynamics, we see that in this case the exit set is reached in a fast manner. The black cross is the initial value.}\label{fig:2d_trajectories} 
\end{figure}

We consider a two dimensional three-hole potential with one being less significant than the others. Note that this potential has already been used in different contextes and is sometimes referred to as M\"uller-Brown potential, see i.e. \cite{Huo1997,Park2003}. In particular we have $\Omega = [-3,3]^2$ and
\begin{equation*}
 pot(x_1, x_2) = 
3 e^{-x_1^2-(x_2-1/3)^2}-3e^{-x_1^2-(x_2-5/3)^2}-5 e^{-(x_1-1)^2-x_2^2}
- 5 e^{-(x_1+1)^2-x_2^2}+ 0.2 x_1^4+0.2 (x_2-1/3)^4.
\end{equation*}
We choose a ball of radius $0.5$ around $x_{target} \approx [-1.048, -0.042]$ as target set $E$ and set $\Xi = \Omega \setminus E$. Note that $x_{target}$ contains a local minimum of $pot$. Both sets are visualized in figure \ref{fig:pot_E}. 
Again, we compare the performance of controllers with different polynomial degree and see that higher order polynomials substantially increase the performance of the controllers. The best performance is achieved by a controller of polynomial degree $16$. However, lower polynomial degree yields 'good' results as well. For the computation of the controllers we set the number of samples to $N = 10 \cdot dof$, where $dof$ is again the number of degrees of freedom of the tensor train. In this two-dimensional example, we set the rank of the tensor train to be maximal for every polynomial degree. Further, we set $M=100$ and $\tau = 0.1$.
In figure \ref{fig:2d_trajectories}, we plot different trajectories of the dynamical system for both, the uncontrolled and controlled system. Note that the controlled trajectories get steered into the set $E$ within the given time frame of $10$, while the uncontrolled dynamics stay in the minimum on the right. Here, we visually see the effect of controlling this dynamical system. Figure \ref{fig:2d_bars} again visualizes the performance of controllers for different polynomial degrees. 

\begin{figure}[htbp]
\centering
       \includegraphics[width=0.45\textwidth]{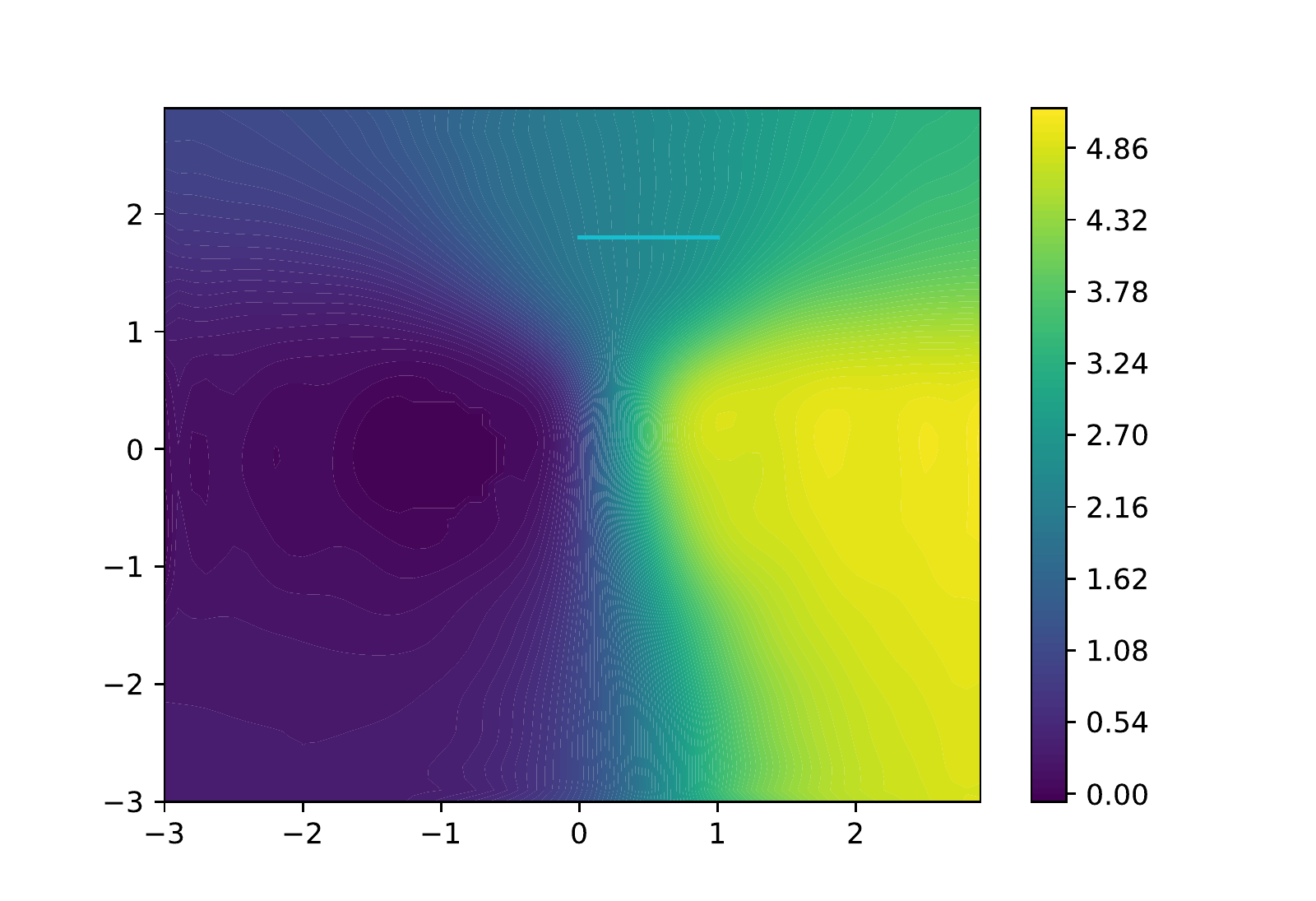}
    \caption{Contour plot of the approximation of the value function. Polynomial degree is $16$.}\label{fig:2d_value_function} 
\end{figure}

As the polynomial degree of $16$ had the best performance, we add contour plots of the value function in Figure \ref{fig:2d_value_function}. In particular compare this Figure to those in \cite{Hartmann2013}, where a similar figure appears in a different context.

\subsection{Test 4: Higher dimensional problem}
\begin{figure}[htbp]
     \subfloat[Red: average distance of the exit set to $1000$ realizations of the uncontrolled dynamics. Black (dots): scatter plot of the distance of $50$ realizations of the uncontrolled dynamics to the exit set.]{%
       \includegraphics[width=0.45\textwidth]{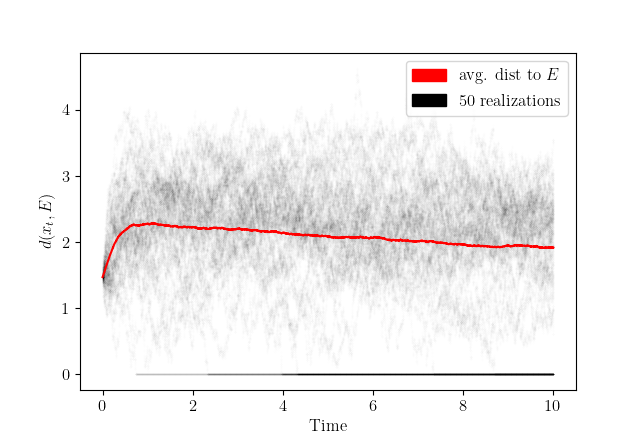}
     }
     \hfill
     \subfloat[Red: average distance of the exit set to $1000$ realizations of the controlled dynamics. Black (dots): scatter plot of the distance of $50$ realizations of the controlled dynamics to the exit set.]{
       \includegraphics[width=0.45\textwidth]{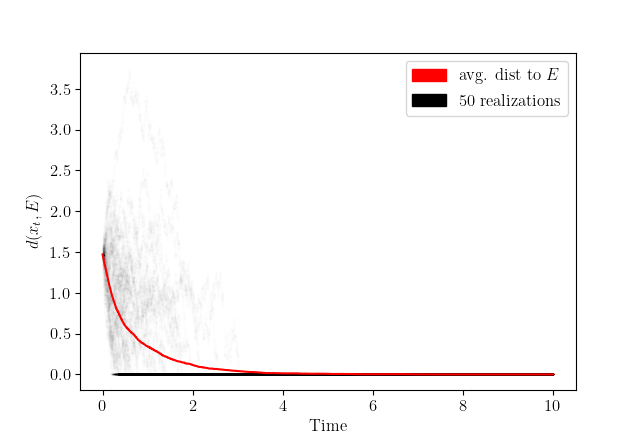}
     }
    \caption{Visualization of the distance of trajectories to the exit set.}\label{fig:ndim_trajectories} 
\end{figure}

We consider the multi-dimensional Double Well potential
\[ pot(x) = \sum_{i=1}^n \kappa_i (x_i^2 - 1)^2. \]
Note that this potential has $2^n$ local minima and the choice of $\kappa$ determines their metastability. In our test we use $\kappa_i = 5$ for all $i$. The exit set $E$ is $B_{0.5 \sqrt{n}}(1)$, the ball of radius $0.5 \sqrt n $ around $[1, \dots, 1] \in \mathbb R^n$. We further choose $\Omega = [- \frac{\pi}{2}, \frac \pi 2 ]^2$. Note that while the radius of the exit set seems large, for the case $n = 6$, its volume is only $1.8\%$ of the volume of $\Omega$. We further stress, that while in the dynamics the dimensions are independent from each other, we cannot expect the value function and thus also the policy to have a strict separation in the dimensions. Here, the curse of dimensions comes into play. Choosing a polynomial degree of $6$ and a tensor train rank of $[5,5,5,5,5]$ allows us to reduce the ansatz space from $46656$ degrees of freedom to $770$. Again, setting $N = 10 \cdot dof$, $M=100$ and $\tau = 0.1$ we visualize the resulting controller in figure \ref{fig:ndim_trajectories}. We start the trajectory at $[-1, \cdots -1] \in \mathbb R^n$. We see that the uncontrolled dynamics approach the exit set slowly, while the controlled dynamics most trajectories have reached the exit set at time $3$. The resulting cost of the controlled dynamics is $4.44$ and the predicted cost $4.74$. Note that computing the average cost for the uncontrolled dynamics is not feasible because of the high metastability of the minima, as seen in Figure \ref{fig:ndim_trajectories}.
\section*{Conclusion and Outlook}
We have considered a stochastic optimal control problem. In the SDE the control $u$ enters as an affine function in $u$ and in its corresponding cost functional $\cal J$ quadratically. 
We solved this problem  by using approximative Policy Iteration whereby we used dynamical programming  
with the linearized Bellman equation resulting in a linear operator equation (Koopman operator). The SDE was discretized by the Euler Mayurama method.
The optimality condition has been derived from the 
HJB equation. For the numerical solution we employed tree based tensor approximations in the subspace of tensor product polynomials. For the computation of the Least Squares risk functional $\cal R$ we used Monte Carlo integration. We have provided successful numerical test for moderate dimensions.

The Least-Squares method allows for incorporating additional penalty terms in (\ref{eq:loss_funct}, \ref{eq:emp-risk}), which might be zero for the exact solution. 
Potential approaches include the following.
\begin{itemize}
    \item The boundary condition $v^\alpha = 0$ on $\partial \Xi$, by sampling the boundary and penalizing $v^\alpha$ at these sample points.
    \item Incorporating the linearized HJB, cf. appendix \ref{sect:HJB} and \eqref{eq:stoch_lin_realhjb} in the risk functional
\end{itemize}
The most time consuming step is the generation of many different paths (trajectories) that are the numerical solution of the SDE.
The above penalty terms yield extra information about the system with low additional computational cost.
In the next future we want to pursue this direction.


\printbibliography

\appendix

\section{Alternate approach via solving the HJB}\label{sect:HJB}
In this section we consider an alternative approach to finding the value function. Indeed, instead of considering the Bellman equation, one can instead consider the Hamilton-Jacobi-Bellman (HJB) equation. It has the form \citep{peng_HJB, Buckdahn2016}
\begin{align}
 \sigma^2(x) \Delta v^*(x) + \min_{\alpha \in F}\{ \nabla v^*( x)\cdot (b(x)+g(x)\alpha(x)) +r(x,\alpha( x))\}&=0 \text{ on } \Xi \label{eq:HJB-stoch} \\
 v^* &= 0 \text{ on } \partial \Xi
\end{align}
with a Dirichlet boundary condition. For the exact value function, 
and the present reward  $ r^\alpha(x) = c(x) + \alpha(x)^T B \alpha(x)$, the minimization within \eqref{eq:HJB-stoch} w.r.t. the parameter $\alpha $ 
can be carried out explicitely. This yields optimality condition for  the optimal policy  (feedback law) given by \cite{DeterministicHJB} 
\begin{equation}\label{eq:HJB-opt}
\alpha^*( x) = - \frac 1 2 B^{-1} g( x)^T \nabla v^*( x).
\end{equation}
Denoting $f^\alpha ( x): = b(x)+g(x)\alpha(x)$ and $r^\alpha( x) = r(x, \alpha( x))$ and $r_{t}^\alpha(x) : = r^\alpha(\Phi_{t}^\alpha(x))$ the corresponding HJB in coupled form is
\begin{align}
    0 &=  \sigma^2(x) \Delta v^*(x)+ \nabla v^*(x) \cdot f^{\alpha^*}(x) + r^{\alpha^*}(x)\label{eq:HJB}\\
    \alpha^*( x) &= - \frac 1 2 B^{-1} g( x)^T \nabla v^*( x) .\notag
\end{align}
Similar to the approach in Section \ref{sect:exit}
for given $\alpha$ we have the linear PDE 
    \begin{equation}\label{eq:stoch_lin_realhjb}
       0 =  \sigma^2 \Delta v^\alpha+ \nabla v^\alpha \cdot f^{\alpha} + r^{\alpha}
    \end{equation}
to compute the policy evaluation function $v^\alpha$. Note that this equation corresponds to \eqref{eq:lin_bellman} in the Bellman setting.

As the Policy Iteration algorithm is based on computing $v^\alpha$, one can exchange \eqref{eq:stoch_lin_hjb} with \eqref{eq:stoch_lin_realhjb} within the Policy Iteration algorithm.
This linearized HJB equation is a stationary inhomogenous backward Kolmogorov equation and thus a deterministic linear PDE in a possibly high-dimensional space. The stochastic nature of the underlying dynamical system is expressed in the additional viscosity term, i.e. Laplacian $\Delta$. 
Indeed, the computational bottleneck is the numerical solution of this high dimensional PDE.
Note that the direct connection between the Koopman operator correspoding to the linearized Bellman equation and the linearized HJB equation is that the backward Kolmogorov operator is the generator of the Koopman operator semi-group, see i.e. \cite{pavliotis2014stochastic}.

\section{Reinterpretation of the finite exit time problem as an infinite horizon problem}
In this section we give an informal way to interpret the finite exit time problem as an non smooth infinite horizon optimal control problem. 
Taking this viewpoint has the advantage that standard HJB theory can be (formally) applied and that the usual Koopman operator without the stop time $\eta$ can be definied. 
Note that our implementation we use this viewpoint. To this end we restate the finite exit time problem
\begin{align*}
    {\mathcal{J}}(x,u) &= \mathbb E \int_{0}^{\eta}  c(X_t) + u_t^T B u_t dt,
\end{align*}
subject to 
\begin{align*}
    d X_t &= b(X_t)dt + \sigma(X_t) d W_t + g(X_t) u_t dt \\
    X_0 &= x
\end{align*}
and 
$$\eta = \inf \{t > 0 | X_t \not\in \Xi \}.$$
Denoting the characteristic function on $\Xi$ by $\chi_\Xi$ we can represent this problem as an infinite horizon problem in the following way.
\begin{align*}
    {\mathcal{J}}(x,u) &= \mathbb E \int_{0}^{\infty}  \chi_\Xi(X_t) (c(X_t) + u_t^T B u_t) dt,
\end{align*}
subject to 
\begin{align*}
    d X_t &= \chi_\Xi(X_t) (b(X_t)dt + \sigma(X_t) d W_t + g(X_t) u_t dt) \\
    X_0 &= x.
\end{align*}
Note that the characteristic function basically sets the running cost and the dynamics to zero once the state is steered out of $\Xi$. Moreover, the dynamics and the cost functional is then non smooth.

\end{document}